\documentclass[leqno,12pt]{amsart}
\usepackage{amsmath,amstext,amssymb,amsopn,amsthm,mathrsfs}

\title[Fractional Integration and Fractional Differentiation]
{ Fractional Integration and Fractional Differentiation for $d$-dimensional
Jacobi Expansions}
\author[C.Balderrama]{Cristina Balderrama$^{a,\ast}$}
\author[W.Urbina]{Wilfredo~O.~Urbina~R{\it $^{a,b}$}.}

\date{}

\theoremstyle{plain}
\newtheorem{thm}{Theorem}[section]

\newtheorem{lem}[thm]{Lemma}
\newtheorem{prop}[thm]{Proposition}

\theoremstyle{definition}

\newcommand{\RR}{\mathbb{R}}

\newcommand{\x}{{\alpha,\beta}}
\newcommand{\dem}{\medskip \par \noindent \mbox{\bf Proof. }} \def\ep{\hfill{$\Box $}}

\begin{document}
\maketitle

\begin{center}
{\it

 $^a$ Departamento de Matem\'{a}ticas, Facultad de Ciencias, UCV.
     Apartado 40009, Los Chaguaramos, Caracas 1041-A Venezuela \\
 $^b$Department of Mathematics and Statistics, University of New Mexico, Albuquerque, NM
    87131, USA \\
  
}
\end{center}

\footnote{$\ast$ Corresponding author\\
{\it E-mail addresses}: cbalde@euler.ciens.ucv.ve (C. Balderrama),
wurbina@euler.ciens.ucv.ve (W. Urbina). }

\begin{abstract}
In this paper we consider  an alternative orthogonal decomposition of the space $L^2$
associated to the $d$-dimensional Jacobi measure in order to obtain an
analogous result to P.A. Meyer's Multipliers Theorem for $d$-dimensional
Jacobi expansions. Then we define and study the Fractional Integral,
the Fractional Derivative and the Bessel potentials induced by the Jacobi
operator. We also obtain a characterization of the Sobolev or potential
spaces and a version of Calder\'on's reproduction formula for the
$d$-dimensional Jacobi measure.
\vspace{10pt} \\
R\'ESUM\'E.
Dans cet article nous consid\'erons une d\'ecomposition
orthogonale alternative de l'espace $L^2$ associ\'ee \`a la mesure
de Jacobi $d$-dimensionelle afin d'obtenir de r\'esultat analogues
au Th\'eor\`eme des Multiplicateurs de P.A. Meyer pour les
d\'eveloppements $d$-dimensionnels de Jacobi. Nous d\'efinissons
et \'etudions l'integral Fractionnaire, la d\'eriv\'ee
Fractionnaire et les potentiels de Bessel induits par l'operateur
de Jacobi. Nous obtenons \`egalement une charact\'erisation des
espaces de Sobolev ou potetiel de Jacobi et une version de la
formule de reproduction de Calder\'on pour la mesure de Jacobi
$d$-dimensionelle.

\end{abstract}

\section{Introduction}

For the parameters $\alpha = ( \alpha_1, \alpha_2, \cdots,
\alpha_d),$ $\beta = ( \beta_1, \beta_2, \cdots, \beta_d)$,  in
$\mathbb{R}^d$ $\alpha_i,\beta_i > -1$ let us consider the
(normalized) Jacobi measure on $[-1,1]^d$ defined as
\begin{eqnarray}
 \mu^d_{\alpha,\beta}(dx) =
\prod_{i=1}^d
\frac{1}{2^{\alpha_i+\beta_i+1}B(\alpha_i+1,\beta_i+1)}
(1-x_i)^{\alpha_i}(1+x_i)^{\beta_i}dx_i.
 \end{eqnarray}
This normalization gives a probability measure. It is not usually
considered in classical orthogonal polynomial theory.

The $d$-dimensional Jacobi operator is given by
\begin{eqnarray}
\mathcal{L}^{\alpha,\beta} =  \sum^{d}_{i=1} \bigg[
(1-x_i^2)\frac{\partial^2}{\partial x_i^2} +  (\beta_i
-\alpha_i-\left(\alpha_i +\beta_i +2\right)x_i)
\frac{\partial}{\partial x_i} \bigg],
\end{eqnarray}
It is not difficult to see that this is formally a symmetric operator in the space
$L^2([-1,1]^d, \mu^d_\x).$

For a multi--index $\kappa=(\kappa_1,\ldots, \kappa_d) \in
\mathbb{N}^d$ let $ {\vec p_\kappa^{\,\,\x}}$ be the normalized Jacobi
polynomial of order $\kappa,$ defined for $x=(x_1, x_2, \cdots, x_d) \in {\mathbb R}^d$ as
$$ {\vec p_\kappa^{\,\,\x}}(x) = \prod_{i=1}^d p_{\kappa_i}^{\alpha_i,\beta_i}(x_i) $$
where $p_n^\x$ for $n\in \mathbb{N}$ and $\x \in \mathbb{R},$
$\x>-1$, is the one dimensional normalized Jacobi polynomial that
can be defined using Rodrigues formula (see \cite{sz}),
$$ (1-x)^\alpha (1+x)^\beta c_n p_n^\x(x)
= {(-1)^n \over 2^n n!} {\partial^n \over \partial x^n }\left\{
(1-x)^{\alpha +n} (1+x)^{\beta+n} \right\}, x\in [-1,1] .$$ As the
one dimensional Jacobi polynomials are orthonormal with respect to
the one dimensional (normalized) Jacobi measure on $[-1,1]$
$$ \mu_{\alpha,\beta}(dx) = \frac{1}{2^{\alpha+\beta+1}B(\alpha+1,\beta+1)}
(1-x)^{\alpha}(1+x)^{\beta}dx,$$ it is immediate that the
normalized Jacobi polynomials $\{ {\vec p_\kappa^{\,\,\x}}\}$ are orthonormal
with respect to the $d$-dimensional Jacobi measure. Moreover the
family $\{ {\vec p_\kappa^{\,\,\x}}\}$ is an orthonormal Hilbert basis of
$L^2([-1,1]^d, \mu^d_\x)$.

It is well known that Jacobi polynomials are eigenfunctions of the
Jacobi operator ${\mathcal L}^{\x}$ with eigenvalue
$-\lambda_\kappa=-\sum_{i=1}^d
\kappa_i(\kappa_i+\alpha_i+\beta_i+1),$ that is,
\begin{eqnarray}\label{OpEnPol}
{\mathcal L}^{\x}  {\vec p_\kappa^{\,\,\x}} = -\lambda_\kappa  {\vec p_\kappa^{\,\,\x}}
\end{eqnarray}

The $d$-dimensional Jacobi semigroup $\{T_t^\x\}_{t\geq 0}$ is the
Markov operator semigroup in $L^2([-1,1]^d, \mu^d_\x)$ associated
to the Markov probability kernel semigroup (see
\cite{B},\cite{ba3} or \cite{ur2}).
\begin{eqnarray*}
P^\x(t, x, dy) = \sum_{\kappa\in \mathbb{N}^d} e^{-\lambda_\kappa
t}
 {\vec p_\kappa^{\,\,\x}}(x) {\vec p_\kappa^{\,\,\x}}(y) \mu^d_\x (dy) =
p^\x_d(t,x,y)\mu_\x^d(dy),
\end{eqnarray*}
that is
\begin{eqnarray*}
T_t^\x f(x) = \int_{[-1,1]^d} f(y) P^\x(t, x, dy)= \int_{[-1,1]^d}
f(y)p^\x_d(t,x,y) \mu_\x^d(dy).
\end{eqnarray*}
Unfortunately, there is not a reasonable explicit representation
of the kernel \\$P^\x(t, x, dy),$
but that  is not needed in what follows. Alternatively the $d$-dimensional Jacobi semigroup
can be defined as the tensorization of one
dimensional Jacobi semigroups (cf. \cite{ur2}).

The $d$-dimensional Jacobi semigroup $\{T_t^\x\}_{t\geq 0}$ is a
Markov diffu\-sion semigroup, conservative, symme\-tric, strongly
continuous on $L^p([-1,1]^d,\mu^d_\x)$ of positive contractions on
$L^p$, with infi\-nite\-simal generator $ {\mathcal L}^{\x}. $ By
(\ref{OpEnPol}) we have
\begin{equation} \label{JacenPol}
T_t^\x  {\vec p_\kappa^{\,\,\x}} = e^{-\lambda_\kappa t}  {\vec p_\kappa^{\,\,\x}}.
\end{equation}

It can be proven that for $\alpha = ( \alpha_1, \alpha_2, \cdots,
\alpha_d),$ $\beta = ( \beta_1, \beta_2, \cdots, \beta_d) \in
\mathbb{R}^d$ with $\alpha_i,\beta_i > -{1\over 2}$,
$\{T_t^\x\}_{t\geq 0}$is not only a contraction on
$L^p([-1,1]^d,\mu^d_\x)$ but it is also an hypercontractive semigroup, that
is to say,  for any initial condition
$1<q(0)<\infty$ there exists an increasing function
$q:\mathbb{R}^+ \to [q(0), \infty),$ such that for every $f$ and
all $t \geq 0,$
\begin{eqnarray*}
\|T_t^\x f\|_{q(t)} \leq \|f\|_{q(0)}.
\end{eqnarray*}

The proof of this fact is not very well known and it is an indirect one. It is based on
the fact that the one dimensional Jacobi operator satisfies a Sobolev
inequality, that can be proved  by checking that it satisfies a curvature-dimension
inequality, this result was obtained by D. Bakry in \cite{ba1}.
Then it can be proved  that this implies a logarithmic Sobolev
inequality for the one dimensional Jacobi operator.
As the logarithmic Sobolev inequality is stable under
tensorization, \cite{a}, we have that the $d$-dimensional Jacobi
operator also satisfies a logarithmic Sobolev inequality and then
using  L. Gross' famous result  \cite{gross}, that asserts the
equivalence between the hypercontractiviy property and the
logarithmic Sobolev inequality, the result is obtained.
All the implications between these
functional inequalities and L. Gross' result can be found in
\cite{a}. A detailed proof of the hypercontractivity propertie for
the Jacobi semigroup can be found in \cite{B}, see also \cite{ur2}.

From
now on we will consider only  the Jacobi semigroup for the parameters
$\alpha = ( \alpha_1, \alpha_2, \cdots, \alpha_d),$ $\beta = (
\beta_1, \beta_2, \cdots, \beta_d) \in \mathbb{R}^d$ with
$\alpha_i,\beta_i > -{1\over 2}$.

 For
$0<\delta<1$ we define the generalized $d$-dimensional
Poisson--Jacobi semigroup of order $\delta$,
$\{P^{\x,\delta}_t\}$, as
\begin{eqnarray}\label{subord}
P^{\x,\delta}_tf(x) = \int^{\infty}_0 T_s^\x
f(x)\mu^{\delta}_t(ds).
\end{eqnarray}
where $\{\mu^{\delta}_t\}$ are the stable measures on $[0,\infty)$
of order $\delta$.$^{(*)}$\footnote{$^{(*)}$ Stable measures on $[0,\infty)$ are
Borel measures on $[0,\infty)$ such that its Laplace transform
verify $\int_0^{\infty} e^{-\lambda s} \mu^{\delta}_t(ds) =
e^{-\lambda^{\delta}t}.$ For $\delta$ fixed, $\{\mu^{\delta}_t\}$
form a semigroup with respect to the convolution operation, see \cite{fell}.} The
generalized $d$-dimensional Poisson--Jacobi semigroup of order
$\delta$ is a strongly continuous semigroup on
$L^p([-1,1]^d,\mu^d_\x)$ with infinitesimal generator $(-{\mathcal
L}^\x)^{\delta}$. Again, by (\ref{OpEnPol}) we have that
\begin{equation}\label{JacGenenPol}
P^{\x,\delta}_t  {\vec p_\kappa^{\,\,\x}}= e^{-\lambda_\kappa^\delta t}
 {\vec p_\kappa^{\,\,\x}}.
\end{equation}

In the particular case $\delta=1/2$, we have the $d$-dimensional
Poisson--Jacobi semigroup. As it is relevant in what follows we
will denote it simply by $P_t^\x= P^{\x,1/2}_t.$ In this case we can
explicitly compute $\mu_t^{1/2},$
$$ \mu_t^{1/2}(ds) = {t \over 2\sqrt{\pi}} e^{-t^2/4s} s^{-3/2} ds $$
and we have  Bochner's  subordination formula,
\begin{equation}\label{Pois-Jac}
P_t^\x f(x) = {1\over\sqrt{\pi}} \int_0^\infty {e^{-u} \over
\sqrt{u}} T_{t^2/4u}^\x f(x)du.
\end{equation}

The paper is organized as follows. In the next section we will
give a decomposition of the space $L^2([-1,1]^d, \mu^d_\x),$ that
we call a modified Wiener-Jacobi decomposition. In section
\ref{sec3}, using this decomposition and the hypercontractivity
property of the $d$-dimensional Jacobi semigruop, we present an
analogous of Meyer's Multiplier Theorem, \cite{me3}, for
$d$-dimensional Jacobi expansions  and define and study, as in the
one dimensional case, \cite{bu},  the fractional derivatives the
fractional integrals, the Bessel potencials for the
$d$-dimensional Jacobi operator, and Jacobi Sobolev spaces
associated to the $d$-dimensional Jacobi measure. Finally we also
study the asymptotic behavior of the $d$-dimensional
Poisson-Jacobi semigroup and, as a consequence, we present a
version of Calderon's reproducing formula.

For others expansions in terms of classical orthogonal polynomials
there have been similar notions. In \cite{lu} it was studied the
fractional derivative for the Gaussian measure, that is, in the
case of Hermite polynomial expansions. In this article they also
obtain characterizations of the Gaussian Sobolev spaces and a
version of Calderon's reproducing formula for the Gaussian
measure.

In \cite{gu} the Laguerre polynomial expansions was studied. In this
article the autors also obtain an analogous of P.A. Meyer's Multiplier
Theorem for Laguerre expansions and introduce fractional
derivatives and fractional integrals in this setting. They also
study different Sobolev spaces  associated to Laguerre expansions
and Riesz-Laguerre transforms.

Thus with this paper we complete the study of these notions for classical
orthogonal polynomials.  In \cite{bu} we have studied the one dimensional case for Jacobi expansions and in the present article we extend this notions to higher dimensions. Contrary to the
Hermite and Laguerre cases, in this case the passage from one
dimension to several dimensions is not straight forward, due to the no
linearity of the eigenvalues of the Jacobi operator. This will be
explained with more detail in the next section.

In order to simplify notation, we will always not explicitly refer
the dependency of the dimension $d$. For instance, we will denote
by $\| \cdot\|_p $ the norm in $L^p([-1,1]^d, \mu_\x^d),$ that is,
without explicitly referring the dependency of the dimension $d$.


\section{A modified Wiener-Jacobi decomposition.}\label{sec2}

Let us consider for each $n \geq 0$, $C_n^\x$ is the closed subspace of $L^2([-1,1]^d,\mu^d_\x)$
generated by the linear combinations of $\{ {\vec
p_\kappa^{\,\,\x}}: |\kappa|=n\},$ where, as usual for a
multi-index $\kappa$, $|\kappa|= \sum_{i=1}^d \kappa_i.$
Since $\{ {\vec p_\kappa^{\,\,\x}}\}$ is an orthonormal basis of
$L^2([-1,1]^d, \mu^d_\x),$ we have the orthogonal decomposition
\begin{equation}\label{Caosdecomp}
 L^2([-1,1]^d,\mu^d_\x)=\bigoplus_{n=0}^\infty C_n^\x.
\end{equation}
 This is the Wiener-Jacobi decomposition of $L^2([-1,1]^d,\mu^d_\x),$
 which is analogous to the Wiener descomposition of $L^2(\RR^d,\gamma_d)$
 is the Gaussian case.

For  $f \in L^2([-1,1]^d,\mu^d_\x),$ its Jacobi
expansion is given by
$$f= \sum_{n=0}^\infty \sum_{|\kappa|=n} \hat{f}(\kappa)  {\vec p_\kappa^{\,\,\x}},$$
with $\hat{f}(\kappa)=\int_{[-1,1]^d} f(y) {\vec
p_\kappa^{\,\,\x}} \mu_\x^d(dy).$ Then we have the following spectral decompositions
\begin{eqnarray*}
{\mathcal L}^{\x} f =  \sum_{n=0}^\infty \sum_{|\kappa|=n} (-\lambda_\kappa) \hat{f}(\kappa)  {\vec p_\kappa^{\,\,\x}}  \\
T_t^\x f =  \sum_{n=0}^\infty \sum_{|\kappa|=n}  e^{-\lambda_\kappa t} \hat{f}(\kappa)  {\vec p_\kappa^{\,\,\x}}\\
P^{\x,\delta}_t f=  \sum_{n=0}^\infty \sum_{|\kappa|=n}  e^{-\lambda_\kappa^\delta t} \hat{f}(\kappa)  {\vec p_\kappa^{\,\,\x}}.
\end{eqnarray*}

 As the eigenvalues $\lambda_\kappa$ of the
$d$-dimensional Jacobi operator  do not depend linearly on $|\kappa|$, we do not
have an expresion of the action of ${\mathcal L}^{\x} $,
$T_t^\x$ or $ P_t^\x$ over $f$, in terms of the orthogonal
projections over the subspaces $C_n^\x$, as in the one
dimensional case (see \cite{bu}) or in the case of Hermite or Laguerre polynomial
$d$-dimensional expansions (see \cite{lu}, \cite{gu}, for example).

For that reason we are going to consider, in the same spirit as
the Wiener-Jacobi decomposition, an alternative decomposition of $
L^2([-1,1]^d, \mu^d_\x)$ in order to obtain expresions of  ${\mathcal L}^{\x} $
$T_t^\x f$ and $ P_t^\x$ in terms of the orthogonal projections.

For fixed $\alpha = ( \alpha_1, \alpha_2, \cdots, \alpha_d),$
$\beta = ( \beta_1, \beta_2, \cdots, \beta_d)$,  in $\mathbb{R}^d$
such that $\alpha_i, \beta_i > {-{1\over 2}}$ let us consider the set,
$$R^\x=\left\{ r\in \mathbb{R}^+
: \mbox{there exists}\,(\kappa_1,\ldots,\kappa_n) \in \mathbb{N}^d,  \mbox{with} \,
r=\sum_{i=1}^d \kappa_i(\kappa_i+\alpha_i+\beta_i +1) \right\}.$$
$R^\x$ is a numerable subset of $\mathbb{R}^+$, thus it can be written as
$R^\x=\{r_n\}_{n=0}^\infty$ with $r_0<r_1<\cdots. $ Let
$$A_n^\x=\left\{ \kappa=(\kappa_1,\ldots,\kappa_d) \in
\mathbb{N}^d: \sum_{i=1}^d \kappa_i(\kappa_i+\alpha_i+\beta_i+1) =
r_n \right\}.$$
Note that $A_0^\x=\{(0,\ldots,0)\}$ and that if $\kappa\in A_n^\x$ then $\lambda_\kappa=\sum_{i=1}^d
\kappa_i(\kappa_i+\alpha_i+\beta_i+1) =r_n.$

Let $G_n^\x$ denote the closed subspace of
$L^2([-1,1]^d,\mu_\x^d)$ generated by the linear combinations of
$\{ {\vec p_\kappa^{\,\,\x}}: \kappa\in A_n^\x \}.$ By the orthogonality of the
Jacobi polynomials with respect to $\mu_\x^d$ and the density of the polynomials, it is not difficult
to see that $\{G^\x_n\}$ is an orthogonal decomposition of
$L^2([-1,1]^d,\mu_\x^d)$, that is
\begin{equation}\label{Caosdecomp2}
 L^2([-1,1]^d,\mu_\x^d)=\bigoplus_{n=0}^\infty G^\x_n.
\end{equation}

We call (\ref{Caosdecomp2}) a modified Wiener--Jacobi
decomposition, compare with (\ref{Caosdecomp}).

Let us denote by $J_n^\x$ the orthogonal projection of
$L^2([-1,1]^d,\mu^d_\x)$ onto $G_n^\x.$ Then, for $f\in
L^2([-1,1]^d,\mu^d_\x)$ its Jacobi expansion now can be written as
\begin{equation}\label{descL2}
f=\sum_{n=0}^\infty J_n^\x f
\end{equation}
where
\begin{equation*}
J_n^\x f = \sum_{\kappa\in A_n^\x} \hat{f}(\kappa)  {\vec p_\kappa^{\,\,\x}}
\end{equation*}
with $\hat{f}(\kappa)= \int_{[-1,1]^d} f(x)  {\vec p_\kappa^{\,\,\x}}(x)
\mu_\x(dx)$
 the Jacobi--Fourier coefficient of $f$ for the multi--index $\kappa$.

By (\ref{OpEnPol}), (\ref{JacenPol}), (\ref{JacGenenPol}) we have
that for $f\in L^2([-1,1],\mu^d_\x)$ with Jacobi expansion
$f=\sum_{n=0}^\infty J_n^\x f$,  the action of ${\mathcal L}^{\x} $,
$T_t^\x$ or $ P_t^\x$ over $f$ can now be expressed as
\begin{eqnarray}
{\mathcal L}^{\x} f = \sum_{n=0}^\infty (-r_n) J_n^\x f, \label{ExpOp} \\
T_t^\x f = \sum_{n=0}^\infty e^{-r_n t} J_n^\x f, \label{EppSem} \\
P^{\x,\delta}_t f= \sum_{n=0}^\infty e^{-r_n^\delta t} J_n^\x f.
\label{ExpPoJac}
\end{eqnarray}
Thus using the modified Wiener-Jacobi decomposition
(\ref{Caosdecomp2}) we are able to obtain expansions of ${\mathcal L}^{\x},$ $ T_t^\x $ and
$P^{\x,\delta}_t$ in terms of the orthogonal projections $J_n^\x$. As we have mentioned before,
this can not be done with the usual Wiener-Jacobi decomposition (\ref{Caosdecomp}).

As a consequence of the hypercontractive property of the $d$-dimensional
Jacobi operator we have that the ortogonal projections
$J_n^\x$ can be extended continuously to $L^p([-1,1]^d,\mu^d_\x)$, more formally

\begin{prop} \label{ContProy}
If $1<p<\infty$ then for every $n\in \mathbb{N}$, $J_n^\x$, restricted to the polynomials ${\mathcal P}$, can be
extended to a continuous operator to $L^p([-1,1]^d,\mu^d_\x),$
that will also be denoted as $J_n^\x$, that is, there exists
$C_{n,p}\in \mathbb{R}^+$ such that
$$ \|J^\x_n f\|_p \leq C_{n,p} \|f\|_p, $$
for $f\in L^p([-1,1]^d, \mu^d_\x).$
\end{prop}
\dem  First remember that the polynomials ${\mathcal P}$ are dense in $L^p([-1,1]^d,\mu^d_\x)$, see
\cite{berg}. Now let us consider $p>2$ and for the initial condition $q(0)=2$,
let $t_0$ be a positive number such that $q(t_0) =p$. Taking $f \in {\mathcal P}$,  then by the hypercontractive property, Parseval's identity
and H\"{o}lder's inequality we obtain,
$$ \|T^\x_{t_0} J^\x_n f\|_p \leq \|J_n^\x f\|_2 \leq \|f\|_2 \leq \|f\|_p. $$
Now, as  $T_{t_0}^\x J_n^\x f =  e^{-t_0 r_n} J_n^\x f$ we get
$$\| J_n^\x f \|_{p} \leq C_{n,p} \|f\|_{p},$$
with $C_{n,p}= e^{t_0 r_n}.$ The general result now follows by density.

 Finally,  for $1<p<2$ the result follows by
duality. \ep


\section{The results} \label{sec3}


Giving  a function $\Phi :\mathbb{N} \to \mathbb{R}$ the
multiplier operator associated to $\Phi$ is defined as
\begin{equation}\label{multop}
T_\Phi f = \sum_{n=0}^\infty \Phi(n) J_n^\x f,
\end{equation}
for $f=\sum_{n=0}^\infty J_n^\x f,\in {\mathcal P}$,  a polynomial.

If $\Phi$ is a bounded function, then by Parseval's identity it is inmediate that
$T_\Phi$ is bounded on $L^2([-1,1]^d,\mu^d_\x).$ In the case of
Hermite expansions, the P.A. Meyer's Multiplier Theorem \cite{me3}
gives conditions over $\Phi$ so that the multiplier $T_\Phi$ can
be extended to a continuous  operator on $L^p$ for $p \neq 2$. In
a previous paper \cite{bu}, we have proven an analogous result for one
dimensional Jacobi expansion. Now we are going to present the analogous result
for $d$-dimensional Jacobi expansions. In order to establish this,
we need some previous results.

First we note that for $n\in \mathbb{N}$, $r_n \geq n$. Then, as a
consequence of the $L^p$ continuity of the projections $J_n^\x$
and of the hypercontractivity of the $d$-dimensional Jacobi
operator we have

\begin{lem}
Let $1<p<\infty.$ Then, for each $m\in \mathbb{N}$ there exists a
constant $C_m$ such that
\begin{eqnarray*}
\|T_t^\x(I-J_0^\x-J_1^\x- \cdots - J_{m-1}^\x) f\|_p \leq C_m
e^{-tm}\|f\|_p.
\end{eqnarray*}
\end{lem}
\dem Let $p>2$ and for the initial condition $q(0)=2$, let $t_0$
be a positive number such that $q(t_0) =p$.

If $t\leq t_0$, since $T_t^\x$ is a contraction, by the $L^p$
continuity of the projections $J_n^\x$,
\begin{eqnarray*}
\|T_t^\x (I-J_0^\x- \cdots - J_{m-1}^\x) f\|_p
&\leq& \| (I-J_0^\x- \cdots - J_{m-1}^\x) f \|_p \\
&\leq& \|f\|_p + \sum_{n=0}^{m-1} \|J_n^\x f\|_p\\
&\leq& (1+ \sum_{n=0}^{m-1} e^{t_0 r_n} )\| f\|_p.
\end{eqnarray*}
But since $e^{t_0 r_n} \leq  e^{t_0 r_m} $ for all $0\leq n\leq
m-1$ and $r_m \geq m$ for all $m \geq 1,$ we get
\begin{eqnarray*}
\|T_t^\x (I-J_0^\x- \cdots - J_{m-1}^\x) f\|_p &\leq& (1+ m
e^{t_0 r_m} )\|f\|_p = C_m e^{-t_0 r_m} \|f\|_p\\
&\leq& C_m e^{-t m} \|f\|_p,
\end{eqnarray*}
with $C_m = (1+ m e^{t_0 r_m} )  e^{t_0 m}.$

Now suppose $t >t_0.$ For $ f=\sum_{n=0}^{\infty} J_n^\x f$, by
the hypercontractive proper\-ty,
\begin{eqnarray*}
\|T_{t_0}^\x T_t^\x (I-J_0^\x- \cdots- J_{m-1}^\x) f\|_p^2
&\leq& \|T_t^\x (I-J_0^\x- \cdots - J_{m-1}^\x) f\|_2^2 \\
 &=& \|T_t^\x (\sum_{n=m}^\infty J_n^\x f )
\|_2^2
= \|\sum_{n=m}^\infty e^{-tr_n} J_n^\x f  \|_2^2\\
 &=&  \sum_{n=m}^\infty e^{-2t r_n} \|J_n^\x f  \|_2^2
\leq \sum_{n=m}^\infty e^{-2tn}\| J_n^\x f \|_2^2,
\end{eqnarray*}
as $r_n \geq n$ for all $n\geq1.$ Then,  as  $m \leq n$,
\begin{eqnarray*}
\sum_{n=m}^\infty e^{-2tn}\| J_n^\x f  \|_2^2 &\leq& e^{-2tm}
\sum_{n=0}^\infty \| J_{n+m}^\x f \|_2^2  \leq
e^{-2tm}\sum_{n=0}^\infty \| J_{n}^\x f \|_2^2
= e^{-2tm} \|f\|_2^2 \\
&\leq& e^{-2tm} \|f\|_p^2.
\end{eqnarray*}
Thus
\begin{eqnarray*}
\|T_{t_0}^\x T_t^\x (I-J_0^\x-J_1^\x- \cdots - J_{m-1}^\x) f\|_p
&\leq& e^{-tm} \|f\|_p,
\end{eqnarray*}
and therefore,
\begin{eqnarray*}
\|T_{t}^\x (I-J_0^\x-\cdots - J_{m-1}^\x) f\|_p
&=& \|T_{t_0}^\x T_{t-t_0}^\x (I-J_0^\x- \cdots  - J_{m-1}^\x) f\|_p \\
&\leq& e^{-(t-t_0)m} \|f\|_p = C_m e^{-tm}\|f\|_p,
\end{eqnarray*}
with $C_m= e^{t_0m}.$ For $1<p<2$ the result follows by duality.
\ep

Using (\ref{subord}) and  Minkowski's integral inequality, it is not difficult to see
an analogous result for the generalized Poisson--Jacobi semigroup,
that is, for $1<p<\infty$ and each $m\in \mathbb{N},$ there exists
$C_m$ such that
\begin{eqnarray} \label{Po-Jacproy}
\|P^{\x,\gamma}_t(I-J_0^\x- J_1^\x- \cdots - J_{m-1}^\x)f\|_p \leq
C_m e^{-tm^{\gamma}} \|f\|_p.
\end{eqnarray}

From the generalized Poisson--Jacobi semigroup let us define a new
family of operators $\{P_{k,\gamma,m}^\x\}_{k \in \mathbb{N}}$ by
the formula
\begin{eqnarray*}
P_{k,\gamma,m}^\x f = {1\over (k-1)!} \int_0^\infty t^{k-1}
P^{\x,\gamma}_t(I-J_0^\x- J_1^\x- \cdots - J_{m-1}^\x)f dt.
\end{eqnarray*}
By the preceding lemma and again by  Minkowski's integral
inequality we have the $L^p$-continuity of $P_{k,\gamma,m}^\x$, for every $m\in
\mathbb{N}$, that is to say,
 for $1<p<\infty$ and  then is a constant $C_m$ such that
\begin{eqnarray}\label{PrevMeyerNuevo}
\| P_{k,\gamma,m}^\x f\|_p \leq {C_m \over m^{\gamma k}} \|f\|_p.
\end{eqnarray}

In particular, if we take $n\geq m$ and $ \kappa \in A_n^\x, $
then
$$P^{\x,\gamma}_t(I-J_0^\x- J_1^\x- \cdots - J_{m-1}^\x)
 {\vec p_\kappa^{\,\,\x}}= e^{-r_n^\gamma t}  {\vec p_\kappa^{\,\,\x}},$$
 and thus, for all
$k \in \mathbb{N}$
\begin{equation*}
 P_{k,\gamma,m}^\x  {\vec p_\kappa^{\,\,\x}} = {1\over r_n^{\gamma k}} {\vec p_\kappa^{\,\,\x}}.
\end{equation*}
Therefore, for $f\in L^2([-1,1]^d, \mu_\x^d)$, $n\geq m$ and $k\in
A_n^\x$
\begin{equation}\label{PotGenEnPol}
P_{k,\gamma,m}^\x J_n^\x f = {1\over r_n^{\gamma k}} J_n^\x f
\end{equation}
and if $n<m$, $k\in A_n^\x$
\begin{equation}\label{PotGenEnPolcero}
P_{k,\gamma,m}^\x J_n^\x f = 0.
\end{equation}

We are ready to establish P.A. Meyer's
Multipliers Theorem for $d$-dimensional Jacobi expansions.

\begin{thm} \label{Meyer} If for some $n_0 \in \mathbb{N}$ and
$0<\gamma<1$
$$\Phi(k) = h\left({1\over r_k^\gamma}\right), \qquad k\geq n_0,$$
with $h$ an analytic function in a neighborhood of zero, then
$T_\Phi$, the multiplier operator associated to $\Phi$, (\ref{multop}), admits a
continuous extension to $L^p([-1,1]^d, \mu^d_\x).$
\end{thm}
\dem Let
\begin{eqnarray*}
T_\Phi f = T_\phi^1 f + T_\Phi^2 f = \sum_{k=0}^{n_0-1} \Phi(k)
J_k^\x f + \sum_{k=n_0}^\infty \Phi(k) J_k^\x f.
\end{eqnarray*}

By Lemma \ref{ContProy} we have that 
\begin{eqnarray*}
\|T_\Phi^1 f\|_p \leq \sum_{k=0}^{n_0-1} |\Phi(k)| \|J_k^\x f\|_p
\leq \left( \sum_{k=0}^{n_0-1} |\Phi(k)| C_k \right) \|f\|_p,
\end{eqnarray*}
that is, $T_\Phi^1$ is $L^p$ continuous. It remains to be seen
that $T_\Phi^2$ is also $L^p$ \mbox{continuous.}

By hypothesis  $h$ can be written as $h(x) = \sum_{n=0}^\infty a_n
x^n, $ for $x$ in a neighborhood of zero, then
\begin{eqnarray*}
T_\Phi^2 f = \sum_{k=n_0}^{\infty} \Phi(k) J_k^\x f =
\sum_{k=n_0}^{\infty} h( {1 \over r_k^\gamma}) J_k^\x f =
\sum_{k=n_0}^{\infty}  \sum_{n=0}^\infty a_n {1\over r_k^{\gamma
n}} J_k^\x f,
\end{eqnarray*}
but by (\ref{PotGenEnPol}) and (\ref{PotGenEnPolcero}), for $k
\geq n_0, \, {1\over r_k^{\gamma n}} J_k^\x f =
P_{n,\gamma,n_0}^\x J_k^\x f$, we have
\begin{eqnarray*}
T_\Phi^2 f  &=&  \sum_{k=n_0}^{\infty} \sum_{n=0}^\infty a_n
P_{n,\gamma,n_0}^\x J_k^\x f =  \sum_{n=0}^\infty a_n
\sum_{k=0}^{\infty} P_{n,\gamma,n_0}^\x J_k^\x f \\
&=& \sum_{n=0}^\infty a_n P_{n,\gamma,n_0}^\x \sum_{k=0}^{\infty}
J_k^\x f  = \sum_{n=0}^\infty a_n P_{n,\gamma,n_0}^\x f.
\end{eqnarray*}
Since  $P_{n,\gamma,n_0}^\x$ is $L^p$ continuous,
(\ref{PrevMeyerNuevo}), we obtain,
\begin{eqnarray*}
\|T_\Phi^2 f\|_p \!\!\!\! &\leq& \!\!\!\! \sum_{n=0}^\infty |a_n| \|P_{n,\gamma,n_0}^\x f\|_p \\
\!\!\!\!&\leq&\!\!\!\! \left(\sum_{n=0}^\infty |a_n| C_{n_0}
{1\over n_0^{\gamma n}} \right) \|f\|_p = C_{n_0}
\left(\sum_{n=0}^\infty |a_n| {1\over n_0^{\gamma n}}
\right)\|f\|_p = C_{n_0}h \left({1\over n_0^{\gamma
}}\right)\|f\|_p.
\end{eqnarray*}

Therefore, $T_\Phi$ is continuous in $L^p([-1,1]^d, \mu_\x^d).$
\ep


As in the classical case of the Laplacian \cite{zy} and in the one dimensional Jacobi case \cite{bu},
for $\gamma >0$ we define the fractional integral of order
$\gamma$, $I_\gamma^\x$, with respect to $d$-dimensional Jacobi operator,  as
\begin{equation}
I_\gamma^\x= (-{\mathcal{L}}^{\x})^{-{\gamma / 2}}.
\end{equation}
$I_\gamma^\x$ is also called Riesz potential of order $\gamma$.

Observe that, since zero is an eigenvalue of ${\mathcal{L}}^\x,$
then $I_\gamma^\x$ is not defined over all
$L^2([-1,1]^d,\mu^d_\x).$ Then consider $\Pi_0=I-J_0^\x$ and denote also by
$I_\gamma^\x$ the operator $(-{\mathcal{L}}^\x)^{-{\gamma / 2}}
\Pi_0.$ Then, this operator is well defined over all
$L^2([-1,1]^d,\mu^d_\x).$ In particular, for the Jacobi polynomial
of order $\kappa$ with $\kappa\in A_n^\x$ we have
\begin{eqnarray}\label{IntenPol}
I_\gamma^\x  {\vec p_\kappa^{\,\,\x}} ={1\over \lambda_\kappa^{\gamma /2}}
 {\vec p_\kappa^{\,\,\x}} ={1\over r_n^{\gamma /2}}  {\vec p_\kappa^{\,\,\x}}.
\end{eqnarray}
Thus, for $f$ a polynomial in  $L^2([-1,1]^d,\mu^d_\x)$ with
Jacobi expansion $\sum_{n=0}^\infty J_n^\x f, $ we have
\begin{eqnarray*}
I_\gamma^\x f= \sum_{n=1}^\infty {1\over r_n^{\gamma /2}} J_n^\x
f.
\end{eqnarray*}

Now, for $\kappa\in A_n^\x$ we have that by the change of
variables $s=\lambda_\kappa^{1/2} t$
\begin{eqnarray*}
{1\over \Gamma(\gamma)} \int_0^\infty t^{\gamma-1} P_t^\x
 {\vec p_\kappa^{\,\,\x}} \, dt &=& {1\over \Gamma(\gamma)} \int_0^\infty
t^{\gamma-1}
e^{-\lambda_\kappa^{1/2} t} \, dt \,  {\vec p_\kappa^{\,\,\x}} \\
&=& {1\over \Gamma(\gamma)} {1\over \lambda_\kappa^{\gamma/2}}
\int_0^\infty s^{\gamma-1} e^{-s} \, ds \,  {\vec p_\kappa^{\,\,\x}} = {1\over
\lambda_\kappa^{\gamma/2}} \,  {\vec p_\kappa^{\,\,\x}},
\end{eqnarray*}
where $P_t^\x$ is the Poisson--Jacobi semigroup.

Therefore for the fractional integral of order $\gamma>0$ we have
the integral representation,
\begin{eqnarray} \label{repInt}
I_\gamma^\x f = {1\over \Gamma(\gamma)} \int_0^\infty t^{\gamma-1}
P_t^\x f dt,
\end{eqnarray}
for $f$ polynomial.

As in the one dimensional case, Meyer's multiplier theorem allows
us to extend $I_\gamma^\x$ as a bounded ope\-ra\-tor on
$L^p([-1,1]^d,\mu^d_\x).$

\begin{thm}\label{ContInt}
The the fractional integral of order $\gamma$, $I_\gamma^\x$
admits a continuous extension, that it will also denoted as denote
$I_\gamma^\x$, to $L^p([-1,1]^d,\mu^d_\x).$
\end{thm}
\dem

If $\gamma/2<1$, then $I_\gamma^\x$ is a multiplier with
associated function
$$\Phi(k)=\frac{1}{r_k^{\gamma/2}}=h\left(\frac{1}{r_k^{\gamma/2}}\right)$$
where $h(z)=z,$ which is analytic in a neighborhood of zero. Then
the  results follows immediately by  Meyer's theorem.

Now, if $\gamma/2\geq1$, let us consider $\beta\in \mathbb{R}$,
$0<\beta<1$ and $\delta=\frac{\gamma}{2\beta}.$ Then
$\delta\beta=\frac{\gamma}{2}$. Let $h(z)=z^\delta$,  which is
analytic in a neighborhood of zero. Then we have
$$h\left(\frac{1}{r_k^\beta}\right)
=\frac{1}{r_k^{\delta\beta}}=\frac{1}{r_k^{\gamma/2}}=\Phi(k).$$
Again the results follows applying Meyer's theorem. \ep


The Bessel potential of order $\gamma>0,$ ${\mathcal J}_\gamma^\x$,
associated to the $d$-dimensional Jacobi operator is defined as
\begin{eqnarray}
{\mathcal J}_\gamma^\x  = (I-{\mathcal{L}}^\x)^{-\gamma /2}.
\end{eqnarray}

For the Jacobi polynomial of order $\kappa$ with $\kappa\in
A_n^\x$ we have
\begin{eqnarray*}
{\mathcal J}_\gamma^\x  {\vec p_\kappa^{\,\,\x}} = {1 \over (1+\lambda_\kappa)^{\gamma/2}}
 {\vec p_\kappa^{\,\,\x}} = {1 \over (1+r_n)^{\gamma/2}}  {\vec p_\kappa^{\,\,\x}} ,
\end{eqnarray*}
and, therefore if $f\in L^2([-1,1]^d,\mu^d_\x)$ polynomial with
expansion $\sum_{n=0}^\infty J_n^\x f$
\begin{eqnarray}
{\mathcal J}_\gamma^\x f = \sum_{n=0}^\infty {1 \over (1+r_n)^{\gamma/2}}
J_n^\x f.
\end{eqnarray}

Again Meyer's theorem allows us to extend Bessel potentials to a
\mbox{continuous} operator on $L^p([-1,1]^d,\mu^d_\x),$

\begin{thm}
The operator ${\mathcal J}_\gamma^\x$ admits a continuous extension, that it
will be also denoted as ${\mathcal J}_\gamma^\x,$ to $L^p([-1,1]^d,
\mu^d_\x).$
\end{thm}
\dem The Bessel Potential of order $\gamma$ is a multiplier
associated to the function
$\Phi(k)=\left(\frac{1}{1+r_k}\right)^{\gamma/2}.$ Let $\beta \in
\mathbb{R}$, $\beta>1$ and
$h(z)=\left(\frac{z^\beta}{z^\beta+1}\right)^{\gamma/2}.$ Then $h$
is an analytic function on a neighborhood of zero and
$$h\left(\frac{1}{r_k^{1/\beta}}\right)
=\left(\frac{1}{1+r_k}\right)^{\gamma/2}=\Phi(k).$$ The results
follows applying Meyer's theorem. \ep


Now, again by analogy to the the classical case of the Laplacian \cite{zy}, we define the fractional derivative of order $\gamma >0$,
$D^\x_\gamma$, with respect to the $d$-dimensional Jacobi operator as
\begin{eqnarray}
 D_\gamma^\x = (-{\mathcal{L}}^\x)^{\gamma/2}.
\end{eqnarray}

For the Jacobi polynomial of order $\kappa$ with $\kappa\in
A_n^\x$ we have
\begin{eqnarray}\label{DerenPol}
D_\gamma^\x  {\vec p_\kappa^{\,\,\x}} = \lambda_\kappa^{\gamma/2}  {\vec p_\kappa^{\,\,\x}} =
r_n^{\gamma/2}  {\vec p_\kappa^{\,\,\x}},
\end{eqnarray}
and therefore, by the density of the polynomials in
$L^p([-1,1]^d,\mu^d_\x)$, $1<p<\infty$, $D_\gamma^\x$ can be
extended to $L^p([-1,1]^d,\mu^d_\x).$

Also for the Jacobi polynomial of order $\kappa$ with $\kappa\in
A_n^\x$ by the change of variables $s=\lambda_\kappa^{1/2} t,$
\begin{eqnarray*}
\int_0^\infty  t^{-\gamma-1} (P_t^\x  {\vec p_\kappa^{\,\,\x}} -  {\vec p_\kappa^{\,\,\x}}) dt
&=& \int_0^\infty  t^{-\gamma-1} (e^{-\lambda_\kappa^{1/2} t}-1) dt   {\vec p_\kappa^{\,\,\x}} \\
&=& \lambda_\kappa^{\gamma/2} \int_0^\infty  s^{-\gamma-1}
(e^{-s}-1) ds \,  {\vec p_\kappa^{\,\,\x}} .
\end{eqnarray*}

Therefore for the fractional derivative  of order $0<\gamma<1$ we also
have a integral representation,
\begin{eqnarray} \label{repDer}
D_\gamma^\x f = {1\over c_\gamma} \int_0^\infty t^{-\gamma-1}
(P_t^\x f-f) dt,
\end{eqnarray}
for $f$ polynomial, where $c_\gamma= \int_0^\infty s^{-\gamma-1}
(e^{-s}-1) ds.$

If $f$ is a polynomial, by (\ref{IntenPol}) and (\ref{DerenPol})
we have,
\begin{eqnarray}
I_\gamma^\x (D_\gamma^\x f) = D_\gamma^\x(I_\gamma^\x f) = \Pi_0
f.
\end{eqnarray}


Let us consider now the Jacobi Sobolev spaces or potential spaces. For $1<p<\infty,$ $L^p_\gamma([-1,1]^d,\mu^d_\x),$ the  Jacobi Sobolev space of order $\gamma>0$,
   is defined as
the completion of the set of polynomials ${\mathcal P}$ with respect to the norm
\begin{eqnarray*}
\|f\|_{p,\gamma} := \|(I-{\mathcal{L}}^\x )^{\gamma/2}f\|_p.
\end{eqnarray*}
That is to say $f\in L^p_\gamma([-1,1]^d,\mu^d_\x)$ if, and only
if, there is a sequence of polynomials $\{f_n\}$ such that
$\lim_{n\to \infty} \|f_n-f\|_{p,\gamma} =0.$

As in the classical case, the Jacobi Sobolev space
$L^p_\gamma([-1,1]^d,\mu^d_\x)$ can also be defined as the image
of $L^p([-1,1]^d,\mu^d_\x)$ under the Bessel Potential
${\mathcal J}_\gamma^\x$, that is,
$$ L^p_\gamma([-1,1]^d,\mu^d_\x) = {\mathcal J}_\gamma^\x L^p([-1,1]^d,\mu^d_\x). $$

The next proposition gives us some inclusion properties among
Jacobi Sobolev spaces,

\begin{prop} \label{IncPot} For the Jacobi Sobolev spaces $L^p_\gamma([-1,1]^d,\mu^d_\x)$, we have
\begin{enumerate}
\item[i)] If $p<q, $ then $L^q_\gamma([-1,1]^d,\mu^d_\x) \subseteq
L^p_\gamma([-1,1]^d,\mu^d_\x)$ for each $\gamma>0.$ \item[ii)]
  If $0<\gamma<\delta,$ then $L^p_\delta([-1,1]^d,\mu^d_\x) \subseteq
L^p_\gamma([-1,1]^d,\mu^d_\x)$ for each $0<p<\infty.$
\end{enumerate}
\end{prop}
\dem

{\em i)} For $\gamma$ fixed, it follows immediately by
H\"{o}lder's inequality.

{\em ii)} Let $f$ be a polynomial and consider
$$ \phi= (I-{\mathcal{L}}^\x)^{\delta/2} f = \sum_{n=0}^\infty (1+r_n)^{\delta/2} J_n^\x f, $$
which is also a polynomial. Then $\phi \in
L^p_{\delta}([-1,1]^d,\mu^d_\x),$ $\|\phi\|_p = \|f\|_{p,\delta}$
and $ {\mathcal J}_{(\gamma-\delta)}^\x \phi =
(I-{\mathcal{L}}^\x)^{(\gamma-\delta)/2} \phi =
(I-{\mathcal{L}}^\x)^{\gamma/2} f, $
by the $L^p$-continuity of
Bessel Potentials,
\begin{eqnarray*}
\|f\|_{p,\gamma} = \|(I-{\mathcal{L}}^\x)^{\gamma/2} f \|_p
=\|{\mathcal J}_{(\gamma -\delta)} \phi \|_p \leq C_{p} \|f\|_{p,\delta}.
\end{eqnarray*}

Now let $f\in L^p_\delta([-1,1]^d,\mu^d_\x)$. Then there exists
$g\in L^p([-1,1]^d,\mu^d_\x)$ such that $f= {\mathcal J}_\delta^\x g $ and a
sequence of polynomials $\{g_n\}$ in $L^p([-1,1]^d,\mu^d_\x)$ such
that $\lim_{n\to \infty}\|g_n-g\|_p=0.$ Set $f_n = {\mathcal J}_\delta^\x
g_n. $ Then $\lim_{n\to \infty} \|f_n-f\|_{p,\delta}=0,$ and
\begin{eqnarray*}
\|f_n-f\|_{p,\gamma}&=&\|(I-{\mathcal{L}}^\x)^{\gamma/2}(f_n-f)\|_p
=\|(I-{\mathcal{L}}^\x)^{\gamma/2}(I-{\mathcal{L}}^\x)^{-\delta/2}(g_n-g)\|_p\\
&=& \|(I-{\mathcal{L}}^\x)^{(\gamma-\delta)/2}(g_n-g)\|_p =
\|{\mathcal J}_{(\gamma-\delta)}(g_n-g)\|_p,
\end{eqnarray*}
by the $L^p$ continuity of Bessel Potentials $\lim_{n\to\infty}
\|f_n-f\|_{p,\gamma} =0. $

Therefore,
\begin{eqnarray*}
\|f\|_{p,\gamma} &\leq& \|f_n-f\|_{p,\gamma} +
\|f_n\|_{p,\gamma}\\
&\leq& \|f_n-f\|_{p,\gamma} + \|f_n\|_{p,\delta},
\end{eqnarray*}
taking limit as $n$ goes to infinity, we obtain the result.\ep

Let us consider the space
$$ L_\gamma([-1,1]^d,\mu^d_\x) = \bigcup_{p>1} L^p_\gamma([-1,1]^d,\mu^d_\x). $$
$L_\gamma([-1,1]^d,\mu^d_\x)$ is the natural domain of
$D_\gamma^\x.$ We define it in this space as follows.

Let $f\in L_\gamma([-1,1]^d,\mu^d_\x),$ then there is $p>1$ such
that $f\in L^p_\gamma([-1,1]^d,\mu^d_\x)$ and a sequence $\{f_n\}$
polynomials such that $ \lim_{n\to \infty}f_n=f$ in $
L^p_{\gamma}([-1,1]^d,\mu^d_\x). $ We define for $f\in
L_\gamma([-1,1]^d,\mu^d_\x) $
\begin{eqnarray*}\label{DerenPot}
D_\gamma^\x f= \lim_{n\to \infty } D_\gamma^\x f_n.
\end{eqnarray*}
The next theorem shows that $D_\gamma^\x$ is well defined and also
inequality (\ref{EqNorPot}) gives us a characterization of the
Sobolev spaces,

\begin{thm}
Let $\gamma>0$ and $1<p,q<\infty.$
\begin{enumerate}
\item[i)] If $\{f_n\}$ is a sequence of polynomials such that
$$\lim_{n\to \infty} f_n = f$$
 in $L^p_{\gamma}([-1,1]^d,\mu^d_\x),$
then
\begin{equation}
\lim_{n\to\infty} D_\gamma^\x f_n \in L^p([-1,1]^d,\mu^d_\x),
\end{equation}
 and the limit does not depend on the choice of  the sequence
$\{f_n\}.$

If $f\in L^p_\gamma([-1,1]^d,\mu^d_\x)\bigcap
L^q_\gamma([-1,1]^d,\mu^d_\x),$ then the limit does not depend on
the choice of $p$ or $q$. Thus $D_\gamma^\x$ is well defined on
$L_\gamma([-1,1]^d,\mu^d_\x)$. \item[ii)] $f\in
L^p_\gamma([-1,1]^d,\mu^d_\x)$ if, and only if, $D_\gamma^\x f\in
L^p([-1,1]^d,\mu^d_\x)$. Moreover, there exists positive constants
$A_{p,\gamma}$ and $B_{p,\gamma}$ such that
\begin{eqnarray} \label{EqNorPot}
B_{p,\gamma} \|f\|_{p,\gamma} \leq \|D_\gamma^\x f\|_p \leq
A_{p,\gamma} \|f\|_{p,\gamma}.
\end{eqnarray}
\end{enumerate}
\end{thm}
\dem

{\em ii)} First, let us note that for $f=\sum_{n=0}^\infty
J_n^\x f$ polynomial,
$$ D_\gamma^\x  {\mathcal J}_\gamma^\x f = \sum_{n=0}^\infty \left({r_n \over 1+r_n }\right)^{\gamma/2} J_n^\x f,$$
that is, the operator $ D_\gamma^\x  {\mathcal J}_\gamma^\x$ is a multiplier
with associated function $\Phi(k) = \left({r_k \over 1+ r_k
}\right)^{\gamma/2} = h({1\over r_k}) $ where $h(z)=\left({ 1
\over z+1}\right)^{\gamma/2},$ and therefore by Meyer's theorem it
is $L^p$-continuos.

Let $f$ be a polynomial and let $\phi$ be a polynomial such that
$f={\mathcal J}_\gamma^\x \phi.$ We have that $\|f\|_{p,\gamma}=\|\phi\|_p$
and by the continuity of the operator $D_\gamma^\x  {\mathcal J}_\gamma^\x$
\begin{eqnarray*}
\|D_\gamma^\x f\|_p =\|D_\gamma^\x  {\mathcal J}_\gamma^\x \phi\|_p \leq
A_{p,\gamma} \|\phi\|_p = A_{p,\gamma} \|f\|_{p,\gamma}.
\end{eqnarray*}

To prove the converse, let us suppose that  $f$ polynomial, then
$D_\gamma^\x f$ is also a polynomial, and therefore $D_\gamma^\x
f\in L^p([-1,1]^d,\mu^d_\x).$ Consider
$$\phi = ( I-{\mathcal{L}}^\x )^{\gamma/2}f
= \sum_{k=0}^\infty (1+r_k)^{\gamma/2} J_k^\x f =
\sum_{k=0}^\infty \left( {1+r_k \over r_k} \right)^{\gamma/2}
J_k^\x(D_\gamma^\x f). $$ The mapping
$$ g=\sum_{k=0}^\infty J_k^\x g \mapsto
\sum_{k=0}^\infty \left( {1+r_k \over r_k} \right)^{\gamma/2}
J_k^\x g  $$ is a multiplier with associated function $ \Phi(k)=
\left( {1+ r_k \over r_k} \right)^{\gamma/2} = h({1 \over r_k}) $
where $h(z)=(z+1)^{\gamma/2},$ so by Meyer's theorem, taking
$g=D_\gamma^\x f$ we have
\begin{eqnarray*}
\|f\|_{p,\gamma} = \|\phi\|_p \leq B_{p,\gamma} \|D_\gamma^\x f\|_p.
\end{eqnarray*}
Thus we get (\ref{EqNorPot}) for polynomials.

For  the general
case, $f\in L^p_\gamma([-1,1]^d,\mu^d_\x)$, then there
exists $g\in L^p([-1,1]^d,\mu^d_\x)$ such that $f={\mathcal J}_\gamma^\x g$
and a sequence $\{g_n\}$ of polynomials such that $\lim_{n\to
\infty}\|g_n -g\|_p=0.$ Let $f_n={\mathcal J}_\gamma^\x g_n,$ then
$\lim_{n\to \infty}\|f_n-f\|_{p,\gamma}=0.$ Then, by the
continuity of the operator $D_\gamma^\x {\mathcal J}_\gamma^\x$ and as
$\lim_{n\to \infty}\|g_n -g\|_p=0$,
$$ \lim_{n\to \infty} \|D_\gamma^\x( f_n-f)\|_p =\lim_{n\to \infty}\|D_\gamma^\x {\mathcal J}_\gamma^\x
(g_n-g)\|_p =0. $$ Then, as
$$ B_{p,\gamma} \|f_n\|_{p,\gamma} \leq \|D_\gamma^\x f_n\|_p \leq
A_{p,\gamma} \|f_n\|_{p,\gamma}, $$ the results follows by taking
the limit as $n \rightarrow \infty$ in this inequality.

{\em i)} Let $\{f_n\}$ be a sequence of polynomials such that
$$\lim_{n\to\infty} f_n = f,$$
 in $L^p_\gamma([-1,1]^d,\mu^d_\x)$.
Then, by (\ref{EqNorPot})
\begin{eqnarray*}
\lim_{n\to \infty} \|D_\gamma^\x f_n\|_p \leq B_{p,\gamma}
\lim_{n\to \infty} \|f_n\|_{p,\gamma} = B_{p,\gamma}
\|f\|_{p,\gamma},
\end{eqnarray*}
hence, $\lim_{n\to\infty} D_\gamma^\x f_n \in
L^p([-1,1]^d,\mu^d_\x).$

Now suppose that $\{q_n\}$ is another sequence of polynomials such
that $\lim_{n\to \infty} q_n= f$ in
$L^p_\gamma([-1,1]^d,\mu^d_\x)$. Then $\lim_{n\to \infty}
f_n-q_n=0.$ By (\ref{EqNorPot})
$$ B_{p,\gamma} \|f_n-q_n\|_{p,\gamma} \leq \|D_\gamma^\x f_n- D_\gamma^\x q_n\|_p \leq A_{p,\gamma} \|f_n-q_n\|_{p,\gamma},$$
and now, taking the limit as $n \rightarrow \infty$
we get that $\lim_{n\to
\infty} D_\gamma^\x f_n = \lim_{n\to \infty} D_\gamma^\x q_n$ in
$L^p([-1,1]^d,\mu^d_\x)$ and therefore the limit does not depends on the
choice of the approximating sequence.

Finally, let us suppose that $f\in
L^p_\gamma([-1,1]^d,\mu^d_\x)\bigcap
L^q_\gamma([-1,1]^d,\mu^d_\x)$ and, without loss of generality, let us assume
that $p\leq q$,  then by Proposition \ref{IncPot}, {\em i)},
$L^q_\gamma([-1,1]^d,\mu^d_\x) \subseteq
L^p_\gamma([-1,1]^d,\mu^d_\x)$ and therefore $f\in
L^q_\gamma([-1,1]^d,\mu^d_\x).$ Now, if $\{f_n\}$ is a sequence
of polynomials such that $\lim_{n\to\infty} f_n = f$ in
$L^q_\gamma([-1,1]^d,\mu^d_\x)$ (hence in
$L^p_\gamma([-1,1]^d,\mu^d_\x)$), we have
$$\lim_{n\to \infty} D_\gamma^\x f_n \in L^q([-1,1]^d,\mu^d_\x)
= L^p([-1,1]^d,\mu^d_\x) \bigcap L^q([-1,1]^d,\mu^d_\x).$$
Therefore the limit does not depends on the choice of $p$ or $q$.
\ep


In what follows we will give an alternative representation of
$D_\gamma^\x$ and $I_\gamma^\x$, but first we present a technical
Lemma, were we study the asymptotic behavior of the $d$-dimensional
 Poisson-Jacobi semigroup $\{P_t^\x\}.$

\begin{lem}
If $f\in C^2([-1,1]^d)$ such that $\int_{[-1,1]^d} f(y)\mu_\x^d(dy) =0$  then
\begin{equation} \label{cota1}
\left| {\partial \over \partial t} P_t^\x f(x) \right| \leq
C_{f,\x,d}(1+|x|) e^{-d_\x^{1/2} t}
\end{equation}
where $d_\x = max\{\alpha_j+\beta_j+2: j=1,\ldots, d\}$ and $|x|$
denotes the usual euclidian norm for $x\in \mathbb{R}^d.$

As a consequence  the Poisson-Jacobi semigroup
$\{P_t^\x\}_{t\geq0}$, has exponential decay on $(C_0^\x)^{\perp}
=\bigoplus_{n=1}^{\infty} C_n^\x$. More explicitly, if $f\in C^2([-1,1]^d)$,
such that  $\int_{[-1,1]^d}f(y)\mu_\x^d(dy)=0$ then
\begin{equation}\label{cota2}
|P_t^\x f(x)| \leq C_{f,\x,d}(1+|x|)e^{-d_\x^{1/2}t}.
\end{equation}
\end{lem}
\dem
First, let us see that $\left|\frac{\partial}{\partial t}
T_t^\x f(x)\right| \leq C_{f,\x,d}(1+|x|) e^{-d_\x t}.$ Since
\begin{eqnarray*}
\frac{\partial}{\partial t} T_t^\x f &=& {\mathcal L}^\x T_t^\x f
= \sum^{d}_{i=1} \bigg[ (1-x_i^2)\frac{\partial^2}{\partial x_i^2}
T_t^\x f \\
&& \quad \quad \quad \quad \quad \quad+ (\beta_i -\alpha_i +1-\left(\alpha_i +\beta_i
+2\right)x_i) \frac{\partial}{\partial x_i} T_t^\x f \bigg]
\end{eqnarray*}

it is sufficient to study $ \frac{\partial}{\partial x_i} T_t^\x f $
and $\frac{\partial^2}{\partial x_i^2} T_t^\x f.$

First note that for the one dimensional Jacobi polinomial $p_n^\x$,
 $n\in \mathbb{N}$ and the one dimensional Jacobi semigroup $T_t^\x$, $\x
 >-{1 \over 2},$ we have
 \begin{eqnarray*}
\frac{\partial}{\partial x}T_t^\x f &=& e^{-(\alpha+\beta+2)t}T_t^{\alpha+1,\beta+1} \left(\frac{\partial f}{\partial x} \right) \\
\frac{\partial^2}{\partial x^2}T_t^\x f &=&
e^{-2(\alpha+\beta+3)t}T_t^{\alpha+2,\beta+2}
\left(\frac{\partial^2 f}{\partial x^2}\right)
\end{eqnarray*}
therefore, as for $\kappa=(\kappa_1,\ldots,\kappa_n) \in
\mathbb{N}^d$, $\x \in \left[-{1 \over 2}, \infty\right)^d $
\begin{eqnarray*}
T_t^\x  {\vec p_\kappa^{\,\,\x}} = \prod_{i=1}^d T_t^{\alpha_i, \beta_i}
p_{\kappa_i}^{\alpha_i,\beta_i},
\end{eqnarray*}
\begin{eqnarray*}
{\partial \over \partial x_j} T_t^\x  {\vec p_\kappa^{\,\,\x}}(x)=
e^{-(\alpha_j+\beta_j+2)t} T_t^{\alpha+e_j, \beta+e_j}
\left({\partial \over \partial x_j} {\vec p_\kappa^{\,\,\x}}(x)
\right)
\end{eqnarray*}
and
\begin{equation*}
{\partial^2 \over \partial x_j^2} T_t^\x  {\vec p_\kappa^{\,\,\x}}(x)=
e^{-2(\alpha_j+\beta_j+3)t} T_t^{\alpha+2e_j, \beta+2e_j}
\left({\partial^2 \over \partial
x_j^2} {\vec p_\kappa^{\,\,\x}}(x) \right)
\end{equation*}
with $e_j \in \mathbb{R}^d$ has one in the $j$-th coordinate and
zero elsewhere. Then for $f \in L^2([-1,1]^d,\mu_\x^d)$
\begin{equation*}
{\partial \over \partial x_j}T_t^\x f(x) =
e^{-(\alpha_j+\beta_j+2)t} T_t^{\alpha+e_j, \beta+e_j}
\left({\partial \over \partial x_j}f(x)\right)
\end{equation*}
and
\begin{equation*}
{\partial^2 \over \partial x_j^2}T_t^\x f(x) =
e^{-2(\alpha_j+\beta_j+3)t} T_t^{\alpha+2e_j, \beta+2e_j}
\left({\partial^2 \over \partial x_j^2}f(x)\right).
\end{equation*}
Hence
\begin{eqnarray*}
\left| {\partial \over \partial t}T_t^\x f(x) \right| &=&
\left|{\mathcal L}^\x T_t^\x f(x) \right| \\ \nonumber &\leq&
\sum^{d}_{j=1} \bigg[|1-x_j^2| e^{-2(\alpha_j+\beta_j+3)t}
T_t^{\alpha+2e_j, \beta+2e_j} \left( \left|{\partial^2 \over
\partial x_j^2}f(x)\right|\right)  \\
\nonumber &+& ( \left|\beta_j -\alpha_j +1 \right|+\left(\alpha_j
+\beta_j +2\right) |x_j|) e^{-(\alpha_j+\beta_j+2)t}
T_t^{\alpha+e_j, \beta+e_j} \left( \left|{\partial \over \partial
x_j}f(x)\right|\right) \bigg].
\end{eqnarray*}
As $f\in C^2([-1,1]^d),$ there exists $C_f$ such that
$\left|{\partial \over
\partial x_j}f(x)\right| \leq C_f$ and $ \left|{\partial^2 \over
\partial x_j^2}f(x)\right| \leq C_f. $ Also, for each $j=1,\ldots,d$
$$ |1-x_j^2| \leq |1-x_j||1+x_j| \leq 1+ |x_j| \leq 1+|x|, $$
$$ e^{-2(\alpha_j+\beta_j+3)t} \leq e^{-(\alpha_j+\beta_j+2)t} $$
$$ |\beta_j -\alpha_j +1 |+(\alpha_j
+\beta_j +2) |x_j| \leq C_{\x}(1+|x_j|) \leq C_\x (1+|x|), $$ and then
\begin{eqnarray*}
\left| {\partial \over \partial t}T_t^\x f(x) \right| &=&
\left|{\mathcal L}^\x T_t^\x f(x) \right| \leq C_{f,\x} (1+|x|)
\sum^{d}_{j=1} e^{-(\alpha_j+\beta_j+2)t}  \\
&\leq& C_{f,\x,d} (1+|x|) e^{-d_\x t} .
\end{eqnarray*}

Now
\begin{equation*}
 {\partial \over \partial t}P_t^\x f(x) = {1\over \sqrt{\pi}}
\int_0^\infty {e^{-u} \over \sqrt{u}} {t\over 2u} {\mathcal L}^\x
T_{t^2/4u}^\x f du,
\end{equation*}
then, by the change of variable $u = d_\x s$
\begin{eqnarray*}
\left| {\partial \over \partial t} P_t^\x f(x) \right| &\leq&
C_{f,\x,d}  (1+|x|) \int_0^\infty e^{-u} {t\over 2 \sqrt{\pi}} u^{-3/2}  e^{-d_\x t^2/4u}  du \\
&=& C_{f,\x,d}  (1+|x|) \int_0^\infty
e^{-d_\x s} {t\over 2 \sqrt{\pi}} s^{-3/2}  e^{- t^2/4s} ds \\
&=& C_{f,\x,d}   (1+|x|) \int_0^\infty e^{-d_\x s} \mu_t^{1/2}
(ds) =C_{f,\x,d}  (1+|x|) e^{-d_\x^{1/2} t}.
\end{eqnarray*}
Since we are assuming that $\int_{[-1,1]^d}f(y)\mu_\x^d(dy)=0$,
\begin{equation} \label{d2}
\lim_{t \to \infty} P_t^\x f(x)=0,
\end{equation}
we get
\begin{eqnarray*}
\left|  P_t^\x f(x) \right| &\leq& \int_t^\infty \left|
\frac{\partial}{\partial s} P_s^\x f(x) \right| ds \leq C_{f,\x,d}
\int_{t}^\infty (1+|x|)e^{-d_\x^{1/2}s} ds\\ \nonumber &=&
C_{f,\x,d} (1+|x|)e^{-d_\x^{1/2}t}.
\end{eqnarray*}
\ep

Now, since $\{P_t^\x\}_{t\geq0}$ is an strongly
continuos semigroup, we have
\begin{equation}
\lim_{t \to 0^{+}}P_t^\x f(x)=f(x) \label{d1}
\end{equation}

Let us write
\begin{eqnarray*}
P^{\x}_t f(x)&=&\int^{\infty}_0 T_s^\x f(x)\mu^{1/2}_t(ds)
=  \int_{[-1,1]^d} [\int^{\infty}_0 p^\x_d(s,x,y)\mu^{1/2}_t(ds)] f(y) \mu_\x^d(dy)\\
&=&  \int_{[-1,1]^d} k^\x_d (t, x,y)f(y) \mu_\x^d (dy),
\end{eqnarray*}
where
\begin{equation}\label{e7}
k^\x_d(t,x,y)=\int^{\infty}_0 p^\x_d(s,x,y)\mu^{1/2}_t(ds).
\end{equation}
and define the operator $Q^\x_t$ as
\begin{equation}\label{e8}
Q^\x_t f(x)=-t\frac{\partial}{\partial t}P_t f(x)=\int_{[-1,1]^d}
q^\x_d(t,x,y)f(y) \mu_\x^d(dy),
\end{equation}
with $ q^\x_d(t,x,y)= -t\frac{\partial}{\partial t}k^\x_d(t,x,y).
$

Now we give the alternative representations for $D_\gamma^\x$ and
$I_\gamma^\x$.

\begin{prop} \label{prop1}
Suppose $f$ differentiable with continuos derivatives up to the
second order  
such that $\int_{[-1,1]^d}f(y)\mu_\x^d (dy)=0$, then we have
\begin{equation} \label{repDer2}
-\gamma  D_\gamma^\x f =\frac {1}{c_{\gamma}}\int_0^{\infty}
t^{-\gamma-1}Q_t^\x f dt, \, 0<\gamma<1,
\end{equation}
\begin{equation} \label{repInt2}
\gamma I_\gamma^\x f =\frac {1}{\Gamma(\gamma)}\int_0^{\infty}
t^{\gamma-1}Q_t^\x f dt, \, \gamma >0.
\end{equation}
\end{prop}
\dem Let us start by proving (\ref{repDer2}). Integrating by parts
in (\ref{repDer}) we have
\begin{eqnarray*}
D_\gamma^\x f &=& {1\over c_\gamma} \lim_{\substack{a \to 0^{+} \\
b \to \infty}} \int_a^b t^{\gamma-1}(P_t^\x f- f)dt \\
&=& {1 \over c_\gamma} \lim_{\substack{a \to 0^{+} \\ b \to
\infty}} \left[ \left. -{P_t^\x f-f \over \gamma t^\gamma }
\right|_a^b + {1\over \gamma} \int_a^b t^{-\gamma } {\partial
\over \partial t} P_t^\x f dt \right] \\
&=& -{1 \over \gamma c_\gamma} \int_0^\infty t^{-\gamma-1} Q_t^\x
fdt,
\end{eqnarray*}
since, by (\ref{cota1}),
\begin{eqnarray*}
\lim_{a\to 0^+} \left| {P_a^\x f(x)-f(x) \over a^\gamma}  \right|
&\leq& \lim_{a \to 0^+} {1 \over a^\gamma} \int_0^a \left|
{\partial \over \partial s} P_s^\x f(x) \right| ds \\ &\leq&
C_{f,\x,d} (1+|x|) \lim_{a\to 0^+} {1- e^{-d_\x^{1/2}a} \over
a^\gamma } =0
\end{eqnarray*}
and by (\ref{d1})
\begin{equation*}
\lim_{b\to \infty} {P_b^\x f -f \over b^\gamma} =0.
\end{equation*}

Let us prove (\ref{repInt2}). Integrating by parts in
(\ref{repInt})
\begin{eqnarray*}
I_\gamma^\x f &=& {1 \over \Gamma(\gamma)} \lim_{\substack{a \to
0^{+} \\ b \to \infty}} \int_a^b t^{\gamma -1} P_t^\x f dt \\
&=&{1\over \Gamma(\gamma)} \lim_{\substack{a \to 0^{+} \\ b \to
\infty}} \left[ \left. { t^{\gamma} \over \gamma } P_t^\x
f\right|_a^b - {1 \over \gamma} \int_a^b t^\gamma {\partial \over
\partial t} P_t^\x f dt \right] \\ &=& {1 \over \gamma
\Gamma(\gamma)} \int_0^\infty t^{\gamma-1} Q_t^\x f dt,
\end{eqnarray*}
since, by (\ref{cota1})
\begin{equation*}
\lim_{b\to \infty} \left| b^\gamma P_t^\x f \right| \leq
C_{f,\x,d} (1+|x|) \lim_{b\to \infty} b^\gamma e^{-d_\x^{1/2}b}=0
\end{equation*}
and
\begin{equation*}
\lim_{a\to 0^+} a^{\gamma} P_a^\x f =0.
\end{equation*}
\ep

This representations for $I_\gamma^\x$ and $D_\gamma^\x$ allows us to
obtain  a version of Calder\'on's reproduction formula for the  $d$-dimensional
Jacobi measure,
\begin{thm}\label{teo1}
i) Suppose $f$ $\in L^1([-1,1]^d,\mu_\x^d)$ such that
$\int_{[-1,1]^d}f(y)\mu_\x^d(dy)=0$, then we have
\begin{equation}\label{cald1}
f =\int_0^{\infty} Q_t^\x f\frac{dt}{t}.
\end{equation}
ii) Suppose $f$ a polynomial such that
$\int_{[-1,1]^d}f(y)\mu_\x^d (dy)=0$, then we have
\begin{equation}\label{c1}
f =C_\gamma \int_0^{\infty} \int_0^{\infty} t^{-\gamma}s^\gamma
Q_t^\x \left( Q_s^\x f \right)\frac{ds}{s}\frac{dt}{t},  \ \ \
0<\gamma<1.
\end{equation}
Also,
\begin{equation}\label{c2}
\int_0^{\infty} \int_0^{\infty} t^{-\gamma}s^\gamma Q_t^\x \left(
Q_s^\x f \right)\frac{ds}{s}\frac{dt}{t}=\int_0^{\infty}
u\frac{\partial^2}{\partial u^2}P_u^\x f du.
\end{equation}
\end{thm}
Formula (\ref{c1}) is the version of Calder\'on's reproduction
formula for the $d$-dimen\-si\-onal Jacobi  measure.
\dem To prove
(\ref{cald1}) note that by (\ref{d2}) and (\ref{d1}) we have,
$$\int_0^{\infty}Q_t^\x f\frac{dt}{t}=\lim_{\substack{a \to 0^{+} \\
b \to \infty}}(-\int_a^b \frac{\partial}{\partial t}P_t^\x f dt)
=\lim_{\substack{a \to 0^{+} \\ b \to \infty}}(-P_t^\x f)
\bigr|^{b}_{a} = f.$$

Let us prove (\ref{c1}). Given $f$ a polynomial such that
$\int_{[-1,1]^d}f(y)\mu_\x^d (dy)=0$, by Proposition \ref{prop1},
we have
\begin{eqnarray} \label{der}
 D_\gamma^\x \left( I_{\gamma}^\x f \right)= - \frac{1}{\gamma c_{\gamma}}\int_0^{\infty}t^{-\gamma-1}
 Q_t^\x \left(I_{\gamma}^\x f \right)dt.
\end{eqnarray}
Now, by (\ref{repInt2}) and the linearity of $Q_t^\x$, we have
\[
Q_t^\x \left(I_{\gamma}^\x f \right)  =\frac{1}{\gamma
\Gamma(\gamma)}\int_0^{\infty}s^{\gamma-1}Q_t^\x(Q_s^\x f)(y)ds.
\]
 Substituting in (\ref{der})
\[
f=D_{\gamma}^\x \left( I_{\gamma}^\x f \right) =C_\gamma
\int_0^{\infty}\int_0^{\infty} t^{-\gamma-1} s^{\gamma-1}Q_t^\x
\left(Q_s^\x f \right) dsdt,
\] with $C_\gamma= -{1\over \gamma^2 c_\gamma \Gamma(\gamma)}.$

In order to prove (\ref{c2}), integrating by parts, and by Proposition
\ref{prop1} we have
\begin{eqnarray*}
\int_0^\infty u \frac{\partial^2}{\partial u^2}P_u^\x f(x)du
&=&\left. u\frac{\partial}{\partial u}P_u^\x
f(x)\right|_{0}^\infty
- \int_0^\infty \frac{\partial}{\partial u}P_u^\x f(x) du \\
&=& - \int_0^\infty \frac{\partial}{\partial u}P_u^\x f(x) du  =\left. -P_u^\x f(x)\right|_0^\infty \\
&=& P_0^\x f(x)= f(x).
\end{eqnarray*}
\ep


\begin{thebibliography}{99}

\pagestyle{myheadings}

\bibitem{a}
An\'e C., Blacher\'e, D., Chaifa\"{\i} D. , Foug\`{e}res  P.,
Gentil, I. Malrieu F. , Roberto C., G. Sheffer G.  {\em Sur les
in\'egalit\'es de Sobolev logarithmiques}. Panoramas et
Synth\`{e}ses 10, Soci\'et\'e Math\'ematique de France. Paris.
2002.

\bibitem{B}
Balderrama, C. {\em Sobre el semigrupo de Jacobi}.
Master Thesis. UCV. 2006.

\bibitem{bu}
Balderrama, C., Urbina, W. {\em Fractional Integration and
Fractional Differentiation for Jacobi Expansions}. Submitted
for publication to
Divulgaciones Mate\'aticas. 2006.

\bibitem{ba1}
Bakry, D. {\em Remarques sur les semi-groupes de Jacobi.}  In
Hommage a P.A. Meyer et J. Neveau. {\bf 236}, Asterique, 1996,
23--40.

\bibitem{ba3}
Bakry, D.  {\em Functional inequalities for Markov semigroups}.
Notes of the CIMPA course in Tata Institute, Bombay, November
2002.

\bibitem{bav}
Bavinck, H. {\em A Special Class of Jacobi Series and Some
Applications}. J. Math. Anal. and Applications. {\bf 37}, 1972,
767--797.

\bibitem{berg} Berg C., Reus C. J. P. {\it Density questions in the classical theory of moments.}
Annales de l'institut Fourier, tome 31, 3 (1981), 99-114.

\bibitem{fell}
Feller, W {\em  An Introduction to Probability Theory and applications} Vol 2. John Wiley \& Sons, Inc.
 New York. 1971.


\bibitem{gu}
Graczyk P., Loeb J.J., L\'opez I.A., Nowak A., Urbina W. {\em
Higher order Riesz Transforms, fractional derivatives and Sobolev
spaces for Laguerre expansions}. Math. Pures Appl. (9), {\bf 84}
2005, no. 3, 375--405.

\bibitem{gross}
L{.} Gross, \emph{Logarithmic Sobolev inequalities and
contractivity properties of semigroups}, Dirichlet forms (Varenna,
1992), Lecture Notes \#1563. Springer, Berlin, 1993, p. 54--88.

\bibitem{lu}
L\'opez I.A., Urbina W. {\em Fractional Differentiation for the
Gaussian Measure and applications}. Bull. Sci. math. {\bf 128},
2004, 587--603.

\bibitem{me3}
Meyer, P. A. {\em Transformations de Riesz pour les lois
Gaussiennes.} Lectures Notes in Math  {\bf 1059}, 1984
Springer-Verlag. 179--193.

\bibitem{sz}
Szeg\"o, G. {\em Orthogonal polynomials.} Colloq. Publ. 23. Amer.
Math. Soc. Providence 1959.

\bibitem{ur}
Urbina, W. {\em An\'alisis Arm\'onico Gaussiano.} Trabajo de
ascenso, Facultad de Ciencias, UCV. 1998.

\bibitem{ur2}
Urbina, W. {\em Semigrupos de Polinomios Cl\'asicos y
Desigualdades Funcionales.} Notas de la Escuela
CIMPA--Unesco--Venezuela. M\'erida (2006).

\bibitem{zy}
 Zygmund, A. {\em Trigonometric  Series.} 2nd. ed. Cambridge Univ. Press. Cambridge (1959).
\end{thebibliography}
\end{document}